\documentclass[reqno]{amsart}
\usepackage{latexsym,amssymb}
\usepackage{upref}
\usepackage{graphicx}
\usepackage[colorlinks=printing,backref]{hyperref}      

\newtheorem{theorem}{Theorem}[section]
\newtheorem{lemma}[theorem]{Lemma}
\newtheorem{corollary}[theorem]{Corollary}
\newtheorem{proposition}[theorem]{Proposition}
\theoremstyle{definition}
\newtheorem{definition}[theorem]{Definition}

\newtheorem{remark}[theorem]{Remark}
\numberwithin{equation}{section}

\begin{document}

\title[A priori estimate and quasiconformal mappings]{A priori estimate of gradient of a solution
to certain differential inequality and quasiconformal mappings}

\author{David Kalaj}
\address{University of Montenegro, Faculty of Natural Sciences and
Mathematics, Cetinjski put b.b. 81000, Podgorica, Montenegro.}
\email{davidk@t-com.me}

\begin{abstract}
 We will prove a global estimate for the gradient of the solution to the {\it Poisson
 differential inequality} $|\Delta u(x)|\le a|\nabla u(x)|^2+b$, $x\in B^{n}$, where
$a,b<\infty$ and $u|_{S^{n-1}}\in C^{1,\alpha}(S^{n-1}, \Bbb R^m)$.
If $m=1$ and $a\le (n+1)/(|u|_\infty4n\sqrt n)$, then $|\nabla u|\,
$ is a priori bounded. This generalizes some similar results due to
E. Heinz (\cite{EH}) and Bernstein (\cite{BS}) for the plane. An
application of these results yields the theorem, which is the main
result of the paper: A quasiconformal mapping of the unit ball onto
a domain with $C^2$ smooth boundary, satisfying the Poisson
differential inequality, is Lipschitz continuous. This extends some
results of the author, Mateljevi\'c and Pavlovi\'c from the complex
plane to the space.
\end{abstract}
\maketitle \tableofcontents
\section{Introduction and statement  of main results}
In the paper $B^n$ denotes the unit ball in $\Bbb R^n$, and
$S^{n-1}$ denotes the unit sphere ($n>2$). We consider the vector
norm $|x|=({\sum_{i=1}^n x_i^2})^{1/2}$ and the matrix norm
$|A|=\max\{|Ax|: |x|=1\}$. Let $\Omega\subset \Bbb R^n$ and
$\Omega'\subset \Bbb R^m$ be open sets and let $u:\Omega\to \Omega'$
be a differentiable mapping. By $\nabla u$ we denote its derivative,
i.e. $$\nabla u=\begin{pmatrix}
 D_1 u_{1} & \dots & D_nu_{1} \\
  \vdots & \dots & \vdots \\
  D_1u_{m} & \dots & D_nu_{m}
\end{pmatrix}.$$ If $n=m$, then the Jacobian of $u$ is defined by
$J_u=\det \nabla u$. The Laplacian of a twice differentiable mapping
is defined by $$\Delta u=\sum_{i=1}^n D_{ii}u.$$

The solution of the equation $\Delta u=g$ (in the sense of
distributions see \cite{Her}) in the unit ball, satisfying the
boundary condition $u|_{S^{n-1}}=f\in L^1(S^{n-1})$, is given by
\begin{equation}\label{green}u(x)=\int_{S^{n-1}}P(x,\eta)f(\eta)d\sigma(\eta)-\int_{B^{n}}G(x,y)g(y)dy,\,
|x|<1.\end{equation} Here
\begin{equation}\label{poisson}P(x,\eta)=\frac{1-|x|^2}{|x-\eta|^n}\end{equation}
is the Poisson kernel and $d\sigma$ is the surface $n-1$ dimensional
measure of the Euclidean sphere satisfying the condition:
$\int_{S^{n-1}}P(x,\eta)d\sigma(\eta)\equiv 1$. The first integral
in (\ref{green}) is called the Poisson integral and is usually
denoted by $P[f](x)$. It is a harmonic mapping. The function
\begin{equation}\label{green1}G(x,y)=c_n\left(\frac{1}{|x-y|^{n-2}}
-\frac{1}{(|1+|x|^2|y|^2-2\left<x,y\right>)^{(n-2)/2}}\right),\end{equation}
where
\begin{equation}\label{cen}c_n=\frac{1}{(n-2)\omega_{n-1}}\end{equation} and
$\omega_{n-1}$ is the measure of $S^{n-1}$, is the Green function of
the unit ball. The Poisson kernel and the Green function are
harmonic in $x$.

\begin{definition}A homeomorfism (continuous mapping) $u:\Omega \to \Omega'$
between two open subsets $\Omega$ and $\Omega'$ of the Euclidean
space $\Bbb R^n$ will be called a $K$ ($K\ge 1$) {\it
quasi-conformal } ({\it quasi-regular}) or shortly a q.c. (q.r.)
mapping if:

(i) $u$ is absolutely continuous function in almost every segment
parallel to some of the coordinate axes and there exist partial
derivatives which are locally $L^n$ integrable functions on
$\Omega$. We will write $u \in ACL^n$.

(ii) $u$ satisfies the condition
\begin{equation}\label{quasin}\frac{|\nabla u(x)|^n}{K}\leq |J_u(x)|\leq
Kl(\nabla u(x))^n\ \ \ \left(\frac{|\nabla u(x)|^n}{K}\leq
J_u(x)\leq Kl(\nabla u(x))^n\right)\end{equation} for almost every
$x$ in $\Omega$ where $l(u'(x)):=\inf\{|\nabla
u(x)\zeta|:|\zeta|=1\}$ and $J_u(x)$ is the Jacobian determinant of
$u$ (see \cite{OS} or \cite{vm}).
\end{definition}
We refer also to the monographs \cite{YUR} and \cite{SRI} for the
basic theory of quasiregular mappings.

Notice that the condition $u \in {ACL}^n$ guarantees the existence
of the first derivative of $u$ almost everywhere (see \cite{OS} ).
Moreover $J_u(x)=\det(\nabla u(x))\neq 0$ for a.e. $x\in \Omega$.
For a continuous mapping $u$, the condition (i) is equivalent to the
fact that $u$ belongs to the Sobolev space
$W^1_{n,\mathrm{loc}}(\Omega)$.

For a function (a mapping) $u$ defined in a domain $\Omega$ we
define $|u| = |u|_\infty =\sup\{|u(x)|: x\in \Omega\}$. We say that
$u\in C^{k,\alpha}(\Omega)$, $0<\alpha\le 1$, $k\in\Bbb N$, if
$$|u|_{l,\alpha}:=\sum_{|\beta|\le l}|D^\beta
u|+\sum_{|\beta|=l}\sup_{x,y\in\Omega}|D^\beta u(x)-D^\beta
u(y)|\cdot|x-y|^{-\alpha}< \infty.$$ It follows that for every
$\alpha\in(0,1]$ and $l\in \Bbb N$
\begin{equation}\label{nninequ}
|u|_l:=\sum_{|\beta|\le l}|D^\beta u|\le |u|_{l,\alpha}.
\end{equation}
We first have that for $u\in C^{1,\alpha}(\Omega)$
\begin{equation}\label{liinequ}
|u(x)-u(y)|\le|u|_{1,\alpha}|x-y| \text{ for every $x,y\in\Omega$,}
\end{equation}
and for real $u\in C^{1,\alpha}(\Omega)$
\begin{equation}\label{liinequ1}
|u^2|_1\le 2|u|_0|u|_{1,\alpha}.
\end{equation}
More generally, for every real differentiable function $\tau$ and
real $u\in C^{1,\alpha}(\Omega)$, we have
\begin{equation}\label{liinequ2}
|\tau(u)|_1\le |\tau'(u)|_0|u|_{1,\alpha}.
\end{equation}

Let $\Omega$ have a $C^{k,\alpha}$ boundary $\partial \Omega$.

The norms on the space $C^{k,\alpha}(\partial \Omega)$ can be
defined as follows.  If $u_0\in C^{k,\alpha}(\partial\Omega)$ then
it has a $C^{k,\alpha}$ extension $u$ to the domain $\overline
\Omega$. The norm in $C^{k,\alpha}(\partial\Omega)$ is defined by:
$$|u_0|_{k,\alpha}:=\inf\{|u|_{\Omega,k,\alpha}:u|_{\partial\Omega}=u_0\}.$$
Equipped with this norm the space $C^{k,\alpha}(\partial\Omega)$
becomes a Banach space. See also \cite[p. 42]{HM} for the definition
of an equivalent norm in $C^{k,\alpha}(\partial\Omega)$, by using
the partition of unity.

One of the starting points of this paper is the following theorem
which was one of the main tools in proving some recent results of
the author and Mateljevic (see \cite{kalmat} and \cite{km}).
\begin{proposition}{\bf (Heinz-Bernstein, see
\cite{BS} and \cite{EH})}. \label{heb} Let $s:\overline{\mathbb{U}}
\to \mathbb{R}$ ($s:\overline{\mathbb{U}} \to \overline{B^m}$) be a
continuous function from the closed unit disc
$\overline{\mathbb{U}}$ into the real line (closed unit ball)
satisfying the conditions:
\begin{enumerate}\item $s$ is $C^2$ on ${\mathbb{U}}$,
\item $s_b(\theta)= s(e^{i\theta})$ is $C^2$ with $K=\max_{\varphi\in[0,2\pi)}|\frac{\partial^2s}{\partial\varphi^2}(e^{i\varphi})|,$ and
\item $ |\Delta s| \leq a |\nabla s|^2+b$ on $\mathbb{U}$ for some
constants $a<\infty$ ($a<1/2$ respectively) and $b<\infty$.
\end{enumerate} Then the function $|\nabla s|\,
$ is bounded on $\mathbb{U}$ by a constant $c(a,b,K)$.
\end{proposition}
The Heinz-Bernstein theorem appeared on 1910 in the Bernstein's
paper \cite{BS} and was reproved by E. Heinz on 1956 in \cite{EH}.
This theorem is important in connection with the Dirichlet problem
for the system
\begin{equation}\label{system}
\Delta u=Q\left(\frac{\partial u}{\partial x_1},\frac{\partial
u}{\partial x_2},u,x_1,x_2\right);
u=(u_1(x_1,x_2),\dots,u_m(x_1,x_2)),
\end{equation}
where $Q=(Q_1,\dots,Q_m)$, and $Q_j$ are quadratic polynomials in
the quantities $\frac{\partial u_i}{\partial x_k}$, $i=1,\dots, m$,
$k=1,2$ with the coefficients depending on $u$ and
$(x_1,x_2)\in\Omega$. An example of the system (\ref{system}) is the
system of differential equations that present a regular surface $S$
with fixed mean curvature $H$ with respect to isothermal parameters
$(x_1,x_2)$. Also it is important in the connection with minimal
surfaces and the Monge-Ampere equation.

We will consider the $n$ dimensional generalization of the system
(\ref{system}). Indeed, we will consider a bit more general
situation. Assume that a twice differentiable mapping
$u=u(x_1,\dots,x_n)$ satisfies the following differential
inequality, which will be the main subject of the paper:

\begin{equation}\label{pdein}
|\Delta u|\le a|\nabla u|^2+b \text{ where $a,b>0$}.
\end{equation}

The inequality (\ref{pdein}) will be called {\it the Poisson
differential inequality}.

Recall that the harmonic mapping equations for $u=(u^1,\dots
u^n):\mathcal N\to \mathcal M$ of the Riemann manifold $\mathcal
N=(B^n,(h_{jk}\ ))$ into a Riemann manifold $\mathcal
M=(\Omega,(g_{jk}\ ))$ ($\Omega\subset \Bbb R^n$) are
\begin{equation}\label{cron} |h|^{-1/2}\sum_{\alpha,\beta=1}^n\partial_\alpha(|h|^{1/2}
h^{\alpha \beta}\partial_\beta u^i) + \sum_{\alpha, \beta, k, \ell =
1}^n\Gamma^i_{k\ell}(u)D_\alpha u^k D_\beta u^\ell, \,\, i =
1,\dots,n,\end{equation} where $\Gamma^i {}_{k\ell}$ are Christoffel
Symbols of the metric tensor $(g_{jk}\ )$ in the target space
$\mathcal M$:
$$    \Gamma^i {}_{k\ell}=\frac{1}{2}g^{im} \left(\frac{\partial
g_{mk}}{\partial x^\ell} +
    \frac{\partial g_{m\ell}}{\partial x^k} - \frac{\partial g_{k\ell}}{\partial x^m} \right) =
     {1 \over 2} g^{im} (g_{mk,\ell} + g_{m\ell,k} - g_{k\ell,m}),$$
the matrix $(g^{jk}\ )$ ( $(h^{jk}\ )$) is an inverse of the metric
tensor $(g_{jk}\ )$ ($(h_{jk}\ )$), and $|h|=\det(h_{jk}\ )$. See
e.g. \cite{jost} for this definition.
\begin{remark}\label{to}
{\it If Christoffel Symbols  of the metric tensor $(g_{jk}\ )$ are
bounded in $\mathcal M$, and if the metric in $\mathcal N$ is
conformal and bounded i.e. if  $h_{jk}(x) = \rho(x)\delta_{jk}$,
such that $\rho$ is bounded in $\mathcal N$ then $u$ satisfies
\eqref{pdein}.}
\end{remark}

 Note that the Poisson differential inequality is related to
the problem
\begin{equation}\label{div}-\mathrm{div}(A(\cdot, u)\nabla u) = f(\cdot,
u, \nabla u),\end{equation} where $x\in B^n(r):=rB^n$, $u(x)\in \Bbb
R^m$ and each $A(x,u)$, for $x\in \Bbb R^n$ and $u\in \Bbb R^m$ is
an endomorphism on $\mathrm{Hom}(\Bbb R^n,\Bbb R^m)$ satisfying
uniformly strongly elliptic and uniformly continuous conditions;
moreover $f$ satisfies the following so called {\it natural growth
condition} (see \cite{gia})
\begin{equation}\label{gia}|f(x,u,p)|\le a(r)|p|^2+b\ \ (p\in \mathrm{Hom}(\Bbb R^n,\Bbb R^m)).\end{equation} The problem of interior and boundary regularity of solutions to \eqref{div} has been treated
by many authors, and among many results, it is proved that, every
solution of \eqref{div}, under some structural conditions and some
conditions on the domains and initial data, has H\"older continuous
extension to the boundary and is $C^{1,\alpha}$ inside. See
\cite{gro} for some recent results on this topic and also \cite{gia}
and \cite{HW} for earlier references.

In this paper we generalize Proposition~\ref{heb} for the space.
Namely we prove the theorems: \vspace{0.5cm}
\\
{\bf Theorem~A.}\,\,{\it Let $s: \overline{{B^n}} \to \mathbb{R}^m$
be a continuous function from the closed unit ball
$\overline{{B^n}}$ into $\mathbb{R}^m$ satisfying the conditions:
\begin{enumerate}\item $s$ is $C^2$ on ${B^n}$,
\item $s:S^{n-1}\to \Bbb R^m$ is $C^{1,\alpha}$ and
\item $ |\Delta s| \leq a |\nabla s|^2+b$ on $B^n$ for some
constants $b $ and $a$ satisfying the condition $ 2a\le|s|_\infty :=
\max\{|s(x)|:x\in \overline{B^n}\}.$
\end{enumerate} Then the function $|\nabla s|\,
$ is bounded on $B^n$. If $a\le C_n/|s|_\infty$ where
\begin{equation}\label{cn}
C_n=\frac{(n+1)\sqrt\pi\Gamma((n + 1)/2)}{8n{\Gamma(n/2+1)}}
\end{equation} then $|\nabla s|\, $ is bounded by a
constant $C(a,b,n,M)$, where $M= |\,s|_{S^{n-1}}\,|_{1,\alpha}$.}
\\
\\
{\bf Theorem~B.}\,\,{\it
 Let $s: \overline{{B^n}} \to
\mathbb{R}$ be a continuous function from the closed unit ball
$\overline{{B^n}}$ into $\Bbb R$ satisfying the conditions:
\begin{enumerate}\item $s$ is $C^2$ on ${B^n}$,
\item $s:S^{n-1}\to \Bbb R$ is $C^{1,\alpha}$ and
\item $ |\Delta s| \leq a |\nabla s|^2+b$ on $B^n$ for some
constants $a$ and $b$.
\end{enumerate} Then the function $|\nabla s|\,
$ is bounded on $B^n$.

If $a\le C_n/|s|_\infty$ then $|\nabla s|$ is bounded by a constant
$ C(a,b,n,|\,u|_{S^{n-1}}\,|_{1,\alpha})$.}

\vspace{0.3cm}

According Remark~\ref{to} it follows that Theorem~A is a partial
extension of \cite[Theorem~4, (ii)]{gh} where it is proved a similar
result, for the family of harmonic mappings $u:\mathcal N\to
\mathcal M$, which map the manifold $\mathcal N$ into a regular ball
$B_r(Q)\subset \mathcal M$.

Theorem~A and Theorem~B roughly speaking, assert that every solution
to Poisson differential inequality has a Lipschitz continuous
extension to the boundary, if the boundary data is $C^{1,\alpha}$.

An application of {Theorem~B} yields the following theorem which is
a generalization of analogous theorems for plane domains due to the
author and Mateljevic (see \cite{ka} and \cite{km}). \vspace{0.5cm}
\\
{\bf Theorem~C.}\,\,{\it Let $u:B^n\to \Omega$ be a twice
differentiable quasiconformal mapping of the unit ball onto the
bounded domain $\Omega$ with $C^2$ boundary, satisfying the Poisson
differential inequality. Then $\nabla u$ is bounded and $u$ is
Lipschitz continuous.} \vspace{0.5cm}

\begin{remark}

 Every $C^2$ mapping satisfies locally the Poisson
differential inequality  and is locally quasiconformal provided that
the Jacobian is non vanishing. Thus the family of mappings
satisfying the conditions of Theorem~C is quite large. According to
Fefferman theorem (\cite{fef}) every biholomorphism of the unit ball
onto a domain with smooth boundary is smooth up to the boundary and
therefore quasiconformal. This fact together with the fact that
every holomorphic mapping is harmonic, it follows that Theorem~C can
be also considered as a partial extension of Fefferman theorem.

\end{remark}

The paper contains this introduction and three other sections. In
section~2 we prove some important lemmas. In section~3, by using
Lemma~\ref{lema8} and  some a priori estimates for Poisson equation
proved in \cite{GH}, we prove Theorem~A and Theorem~B, which are a
priori estimates and global estimates for the the solution to
Poisson differential inequality on the space, which are
generalization of analogous classic Heinz-Bernstein theorem for the
plane. In the final section, we previously show that, if $u:B^n\to
\Omega$ is q.c. and satisfies Pisson differential inequality, then
the function $\chi(x)=-d(u(x))$, where
$d(u)=\mathrm{dist}(u,\partial \Omega)$, satisfies as well Pisson
differential inequality in some neighborhood of the boundary. By
using this fact and Theorem~B we prove Theorem~C. This extends some
results of the author, Mateljevic and Pavlovic (\cite{ka},
\cite{trans}, \cite{km}, \cite{kalmat} and \cite{MP}) from the plane
to the space. It is important to notice that, the conformal mappings
and decomposition of planar harmonic mappings as the sum of an
analytic and an anti-analytic function played important role in
establishing some regularity boundary behaviors of q.c. harmonic
mappings in the plane (\cite{MP} and \cite{ka}). This cannot be done
for harmonic mappings defined in the space (see \cite{akm},
\cite{kan} and \cite{tt} for some results on the topic of Euclidean
and hyperbolic q.c. harmonic mappings in the space). The theorem
presented here (Theorem B) made it possible to work on the problem
of q.c. harmonic mappings on the space without employing the
conformal and analytic functions. Notice that, the family of
conformal mappings on the space coincides with M\"obius
transformations. Therefore this family is "smaller" than the family
of conformal mappings in the plane.

\section{Some lemmas}
 We will follow the approach used in \cite{EH}. The following lemma will be an essential tool in proving the main
results of Section~3. It depends upon three lemmas. It is an
extension of a similar result for complex plane treated in
\cite{EH}.
\begin{lemma}[The main lemma]\label{lema8}
Assumptions:

{\bf A1} The mapping $u:D\to {\overline B^m}$ is $C^2$, $D\subset
B^n$ satisfies for $x\in D$ the differential inequality
\begin{equation}\label{212}
|\Delta u|\le a|\nabla u|^2+b
\end{equation}
where $0<a, b<\infty$.

{\bf A2} There exists a real valued function $G(u)$ of class $C^2$
for $|u|\le 1+\epsilon$ ($\epsilon>0$) such that
\begin{equation}\label{214}
|\nabla G|\le \alpha
\end{equation} and the function $\phi(x)=G(u(x))$ satisfies the
differential inequality
\begin{equation}\label{215}
\Delta \phi\ge \beta(|\nabla u|^2)-\gamma
\end{equation} where $\alpha$, $\beta$ and $\gamma$ are positive
constants.

Conclusions:

{\bf C} There exists a fixed positive number
$c'_1=c_1(a,b,\alpha,\beta,\gamma,n,u)$ such that for $x_0\in D$ and
$r_0 = \mathrm{dist}(x_0,\partial D)$, the following inequality
holds:
\begin{equation}\label{216}
|\nabla u(x_0)|\le c'_1\left(1+\frac{\max_{|x-x_0|\le
r_0}|u(x)-u(x_0)|}{r_0}\right).
\end{equation}
If $a$ is small enough ($a$ satisfies the inequality $a\le C_n$ i.e.
(\ref{ieq})) then $c_1'$ can be chosen independently of $u$.
\end{lemma}
\begin{remark}
After writing this paper, the author learned that a similar result
has been obtained in the paper \cite{jk}, for the class of harmonic
mappings. By using the fact that the family of bounded harmonic
mappings is H\"older continuous in compacts for some H\"older
exponent $\sigma<1$ (a result obtained in \cite{hsw}), Jost and
Karcher proved that $c_1'$ can be chosen independently of $u$
without the previous restriction (\cite[Theorem~3.1]{jk}). However,
it seems that our results are new for the class of solutions to
Poisson differential inequality and the author believes that $c_1'$
can be chosen independently of $u$ without the restriction $a\le C_n
=\frac{(n+1)\sqrt\pi\Gamma((n + 1)/2)}{8n{ \Gamma(n/2 +1)}}$.
However to prove the main result (Theorem~C) we only need an
estimate that is not necessarily a priori (see Theorem~B).
\end{remark}

We will prove the lemma~\ref{lema8} by using the following three
lemmas:

\begin{lemma}\label{lema8'}
Let $u$ satisfy the hypotheses of the lemma~\ref{lema8} and let the
ball $|y-x|\le \rho$ be contained in $D$. Then we have for
$0<\rho_1<\rho$ the inequality
\begin{equation}\label{217}
\begin{split}
\int_{|y-x|\le \rho_1}{c_n}&|\nabla u|^2dy
\\&\le \frac{{\rho_1}^{n-2}\rho^{n-2}}{\rho^{n-2}-\rho_1^{n-2}}
\left(\frac{\gamma
\rho^2}{2n\beta}+\frac{\alpha}{\beta}\max_{|y-x|=\rho}|u(y)-u(x)|\right).
\end{split}
\end{equation}

\end{lemma}

\begin{proof}
By using (\ref{green}) we obtain for $|x|<1-\rho$
\begin{equation}\label{218}\begin{split}
\int_{S^{n-1}}&(u(x+\rho\eta)-u(x))d\sigma(\eta)\\&=\int_{|y-x|\le
\rho}\left(\frac{c_n}{|y-x|^{n-2}}-\frac{c_n}{\rho^{n-2}}\right)g(y)dy.
\end{split}
\end{equation}
If we now apply the identity (\ref{218}) to the mapping
$\phi(x)=G(u(x))$, by using $|\nabla \phi|\le a$, we obtain the
inequality:
\begin{equation}\label{219}\begin{split}
\int_{|y-x|\le
\rho}\left(\frac{c_n}{|y-x|^{n-2}}-\frac{c_n}{\rho^{n-2}}\right)\Delta
\phi dy&\\\le \alpha
\int_{S^{n-1}}|u(x+\rho\eta)-u(x)|d\sigma(\eta)&\le \alpha
\max_{|x-y|=\rho}|u(y)-u(x)|.
\end{split}
\end{equation}
On the other hand from (\ref{215}) we deduce
\begin{equation}\label{2110}\begin{split}
\int_{|y-x|\le
\rho}\Big(\frac{c_n}{|y-x|^{n-2}}&-\frac{c_n}{\rho^{n-2}}\Big)\Delta
\phi dy\\ &\ge \beta\int_{|y-x|\le
\rho}\left(\frac{c_n}{|y-x|^{n-2}}-\frac{c_n}{\rho^{n-2}}\right)|\nabla
u|^2dy\\& -\gamma \int_{|y-x|\le
\rho}\left(\frac{c_n}{|y-x|^{n-2}}-\frac{c_n}{\rho^{n-2}}\right)dy\\&
\ge \beta\int_{|y-x|\le
\rho}\left(\frac{c_n}{|y-x|^{n-2}}-\frac{c_n}{\rho^{n-2}}\right)|\nabla
u|^2dy\\&-\frac{\gamma\rho^2}{2n}.
\end{split}
\end{equation}
Combining this inequality with (\ref{219}) we obtain
\begin{equation}\label{2111}\begin{split}
\int_{|y-x|\le
\rho}\left(\frac{c_n}{|y-x|^{n-2}}-\frac{c_n}{\rho^{n-2}}\right)&|\nabla
u|^2dy\\ &\le\frac{\gamma\rho^2}{2\beta
n}+\frac{\alpha}{\beta}\max_{|y-x|=\rho}|u(y)-u(x)|.
\end{split}
\end{equation}
Now let $0<\rho_1<\rho$. From (\ref{2111}) we get
\[\begin{split}\left(\frac{1}{\rho_1^{n-2}}-\frac{1}{\rho^{n-2}}\right)&\int_{|y-x|\le \rho_1}{c_n}|\nabla
u|^2dy\\ &\le\int_{|y-x|\le
\rho}\left(\frac{c_n}{|y-x|^{n-2}}-\frac{c_n}{\rho^{n-2}}\right)|\nabla
u|^2dy\\&\le \frac{\gamma\rho^2}{2\beta
n}+\frac{\alpha}{\beta}\max_{|y-x|=\rho}|u(y)-u(x)|
\end{split}\]
and therefore $$ \int_{|y-x|\le \rho_1}{c_n}|\nabla u|^2dy\le
\frac{{\rho_1}^{n-2}\rho^{n-2}}{\rho^{n-2}-\rho_1^{n-2}}
\left(\frac{\gamma\rho^2}{2\beta
n}+\frac{\alpha}{\beta}\max_{|y-x|=\rho}|u(y)-u(x)|\right).$$
\end{proof}
\begin{lemma}\label{lemma8''}
Let $Y: D\to B^m$ be a $C^2$ mapping of a domain $D\subset B^n$. Let
$B^n(x_0,\rho)\subset D$ and let $Z\in \Bbb R^m$ be any constant
vector ($n\ge 3$, $m\in \Bbb N$). Then we have the estimate:

\begin{equation}\label{2112}\begin{split}
|\nabla Y(x_0)|&\le
\frac{n}{\rho^n}\int_{|y-x_0|=\rho}|Y(y)-Z|d\sigma(y)\\&+\frac{1}{\omega_{n-1}}\int_{|y-x_0|\le\rho}\left(\frac{1}{|y-x_0|^{n-1}}
-\frac{|y-x_0|}{\rho^n}\right) |\Delta Y| dy,
\end{split}
\end{equation}

and

\begin{equation}\label{shtune}\begin{split}
|\nabla Y(x_0)|&\le
\frac{\gamma_n}{\rho}+\frac{1}{\omega_{n-1}}\int_{|y-x_0|\le\rho}\left(\frac{1}{|y-x_0|^{n-1}}
-\frac{|y-x_0|}{\rho^n}\right) |\Delta Y| dy,
\end{split}
\end{equation}
where $\gamma_n$ is defined in \eqref{nate}.
\end{lemma}
\begin{proof}

Assume that $v\in C^2(\overline{B^n})$.  From
\begin{equation}\label{vh}v(x)=H(x)+K(x):=\int_{S^{n-1}}P(x,\eta)v(\eta)d\sigma(\eta)-\int_{B^{n}}G(x,y)\Delta
v(y)dy\end{equation} where $H$ is a harmonic function, it follows
that
\begin{equation}\label{d3}v'(x)h=\int_{S^{n-1}}P_x(x,\eta)h\cdot
v(\eta)d\sigma(\eta)-\int_{B^{n}}G_x(x,y)h \cdot\Delta
v(y)dy.\end{equation} By differentiating (\ref{poisson}) and
(\ref{green1}) we obtain
$$P_x(x,\eta)h=\frac{-2\left<x,h\right>}{|x-\eta|^n}-\frac{n(1-|x|^2)\left<x-\eta,h\right>}{|x-\eta|^{n+2}}$$
and

$$G_x(x,y)h=c_n\frac{(n-2)\left<x-y,h\right>}{|x-y|^n}-c_n\frac{(n-2)|y|^2\left<x,h\right>-(n-2)\left<y,h\right>}
{(|1+|x|^2|y|^2-2\left<x,y\right>)^{n/2}}.$$ Hence
$$P_x(0,\eta)h=\frac{n\left<\eta,h\right>}{|\eta|^{n+2}}=n\left<\eta,h\right>$$
and
$$G_x(0,y)h=-\frac{1}{\omega_{n-1}}\frac{\left<y,h\right>}{|y|^{n}}+\frac{1}{\omega_{n-1}}\left<y,h\right>.$$
Therefore
\begin{equation}\label{d1}|P_x(0,\eta)h|\le|P_x(0,\eta)|h|=n|h|\end{equation} and
\begin{equation}\label{d2}|G_x(0,y)h|\le|G_x(0,y)|h|=\frac{1}{\omega_{n-1}}({|y|^{1-n}}-|y|)|h|.\end{equation}
By using (\ref{d3}), (\ref{d1}) and (\ref{d2}) we obtain $$|\nabla
v(0)h|\le \int_{S^{n-1}}|P_x(0,\eta)||h|
|v(\eta)|d\sigma(\eta)+\int_{B^{n}}|G_x(0,y)||h||\Delta v(y)|dy.$$
Hence we have
\begin{equation}\label{in}\begin{split}|\nabla v(0)|&\le
n\int_{S^{n-1}}
|v(\eta)|d\sigma(\eta)\\&+\frac{1}{\omega_{n-1}}\int_{B^{n}}({|y|^{1-n}}-|y|)
|\Delta v(y)|dy.\end{split}\end{equation} Let $v(x)=Y(x_0+\rho
x)-Z$. Then $v(0)=Y(x_0)-Z$, $\nabla v(0)=\rho \nabla Y(x_0)$ and
$\Delta v(y)=\rho^2\Delta Y(x_0+\rho y)$. Inserting this into
(\ref{in}) we obtain
\begin{equation}\label{in1}\begin{split}\rho|\nabla Y(x_0)|&\le
|\nabla G(0)|+|\nabla K(0)|\le n\int_{S^{n-1}}
|Y(x_0+\rho\eta)-Z|d\sigma(\eta)\\&+\rho^2\frac{1}{\omega_{n-1}}\int_{B^{n}}({|y|^{1-n}}-|y|)
|\Delta Y(x_0+\rho y)|dy.\end{split}
\end{equation}
Introducing the change of variables $\zeta=x_0+\rho\eta$ in the
first integral and $w=x_0+\rho y$ in the second integral of
(\ref{in1}) we obtain
\begin{equation}\label{ini2}\begin{split} |\nabla Y(x_0)|&\le
\frac{n}{\rho^n}\int_{|\zeta-x_0|=\rho}|Y(\zeta)-Z|d\sigma(\zeta)\\&+\frac{1}{\omega_{n-1}}\int_{|w-x_0|\le\rho}\left(\frac{1}{|w-x_0|^{n-1}}-\frac{|w-x_0|}{\rho^n}\right)
|\Delta Y| dw
\end{split}
\end{equation}
which is identical with (\ref{2112}). To get \eqref{shtune} we do as
follows. We again start by \eqref{vh}. Let $v(x)=Y(x_0+\rho x)$.
Then $H$ is defined by
$$H(x)=\int_{S^{n-1}}P(x,\eta)Y(x_0+\rho\eta)d\sigma(\eta).$$ Applying
the Schwartz lemma (see \cite[Theorem~6.26]{ABR}) to the harmonic
function $H_h:x\mapsto \left<H(x),h\right>$, where $h$ is a unit
vector in $\Bbb R^m$, we obtain that
\begin{equation}\label{nate}\begin{split}|\nabla H(0)|&=\max_{k\in S^{n-1},h\in S^{m-1}}|\left<\nabla H(0)k,h\right>| = \max_{|h|=1}|\nabla H_h(0)|\\&\le \gamma_n:=\frac{2(n-1)\omega_{n-2}}{n\omega_{n-1}} =
\frac{2 \Gamma(1 + n/2)}{\sqrt\pi \Gamma((n + 1)/2)}<\sqrt
n.\end{split}\end{equation} Since, $$Y(x_0+\rho x)  = H(x) + K(x)$$
and $\nabla v(0) = \rho\nabla Y(x_0)$, by using inequality
\eqref{nate} together with the previous estimate for $|\nabla
K(0)|$, it follows \eqref{shtune}.
\end{proof}
\begin{lemma}\label{lem}
Let $u:\overline{B^n}\to\Bbb R^m$ be a continuous mapping. Then
there exists a positive function $\delta_u=\delta_u(\varepsilon)$,
$\varepsilon\in (0,2)$, such that if $x,y\in \overline{B^n}$, and
$|x-y|<\delta_u(\varepsilon)$ then $|u(x)-u(y)|\le \varepsilon$.
\end{lemma}
In \cite[Theorem~3]{hsw} is proved that the family of harmonic
mappings $u:\mathcal N\to \mathcal M$, which map the monifold
$\mathcal N$ into the regular ball $B_r(Q)\subset \mathcal M$ is
uniformly continuous in compact subsets of $\mathcal N$. This
implies that the function $\delta_u(\varepsilon)$ can be chosen
independently on $u$. This fact can improve the conclusion of
Lemma~\ref{lema8} as it is done in \cite[Theorem~3.1]{jk}.

\begin{proof}[Proof of
Lemma~\ref{lema8}] In order to estimate the function $|\nabla u|^2 $
in the ball $|x-x_0|<r_0$ we introduce the quantity
\begin{equation}\label{2114}
M=\max_{|x-x_0|<r_0}(r_0-|x-x_0|)|\nabla u(x)|.
\end{equation}
Obviously there exists a point $x_1:|x_1-x_0|<r_0$ such that
\begin{equation}\label{2115}
M=(r_0-|x_1-x_0|)|\nabla u(x_1)|.
\end{equation}
Let $d=r_0-|x_1-x_0|$ and $\theta\in (0,1)$. If we apply
Lemma~\ref{lemma8''} to the case where $Y(x)=u(x)$ and $Z=u(x_1)$,
$x=x_1$ and $\rho=d\theta$, and use (\ref{2115}), we obtain
\[\begin{split}
\frac{M}{d}&\le\min
\{\frac{\gamma_n}{d\theta},\frac{n}{d^n\theta^n}\int_{|y-x_1|=d\theta}|u(y)-u(x_1)|d\sigma(y)\}\\&
+\frac{1}{\omega_{n-1}}\int_{|y-x_1|\le
d\theta}\left(\frac{1}{|y-x_1|^{n-1}}-\frac{|y-x_1|}{d^n\theta^n}\right)
|\Delta u| dy.
\end{split}
\]
Using now (\ref{212}) we obtain
\begin{equation}\label{2116}\begin{split}
\frac{M}{d}&\le\min
\{\frac{\gamma_n}{d\theta},\frac{n}{d^n\theta^n}\int_{|y-x_1|=d\theta}|u(y)-u(x_1)|d\sigma(y)\}\\&
+(n-2)a c_n\int_{|y-x_1|\le
d\theta}\left(\frac{1}{|y-x_1|^{n-1}}-\frac{|y-x_1|}{d^n\theta^n}\right)|\nabla
u|^2 dy\\&+\frac{b}{\omega_{n-1}} \int_{|y-x_1|\le
d\theta}\left(\frac{1}{|y-x_1|^{n-1}}-\frac{|y-x_1|}{d^n\theta^n}\right)dy.
\end{split}
\end{equation}
We shall now estimate the right hand side of (\ref{2116}). First of
all, according to Lemma~\ref{lem}, we have for every $\varepsilon>0$
and $d\theta<\delta_u(\varepsilon)$ the inequality:
\begin{equation}\label{2117}\displaystyle
\frac{n}{d^n\theta^n}\int_{|y-x_1|=d\theta}|u(y)-u(x_1)|d\sigma(y)\le
\frac{n\varepsilon}{d\theta}.
\end{equation}
 On the other hand
\begin{equation}\label{2118}
\frac{b}{\omega_{n-1}}\int_{|y-x_1|\le
d\theta}\left(\frac{1}{|y-x_1|^{n-1}}-\frac{|y-x_1|}{d^n\theta^n}\right)dy=b\frac{n}{n+1}d\theta.
\end{equation}
Next let $\lambda$ be a real number such that $0<\lambda<\theta$.
Then  we have the inequality
\begin{equation}\label{2119}\begin{split}
&\frac{a}{\omega_{n-1}}\int_{|y-x_1|\le
d\theta}\left(\frac{1}{|y-x_1|^{n-1}}-\frac{|y-x_1|}{d^n\theta^n}\right)|\nabla
u|^2 dy\\&\le\frac{a}{\omega_{n-1}}\int_{|y-x_1|\le
d\lambda}\left(\frac{1}{|y-x_1|^{n-1}}-\frac{|y-x_1|}{d^n\theta^n}\right)|\nabla
u|^2 dy
\\&+(n-2)ad\lambda c_n\left(\frac{1}{d^{n}\lambda^{n}}-\frac{1}{d^{n}\theta^{n}}\right)\int_{|y-x_1|\le d\theta}{|\nabla
u|^2}dy.
\end{split}
\end{equation}
In order to estimate the right hand side of this inequality we first
observe that, on account of (\ref{2115}) we have for $|x-x_1|\le
d\lambda$ the estimate $$|\nabla u|^2\le
\frac{M^2}{d^2(1-\lambda)^2},$$ and therefore
\begin{equation}\label{2120}\begin{split}
\frac{a}{\omega_{n-1}}&\int_{|y-x_1|\le
d\lambda}\left(\frac{1}{|y-x_1|^{n-1}}-\frac{|y-x_1|}{d^n\theta^n}\right)|\nabla
u|^2dy\\&\le
\frac{a}{\omega_{n-1}}\frac{M^2}{d^2(1-\lambda)^2}(\lambda
d-\frac{\lambda^{n+1}d}{(n+1)\theta^n})\omega_{n-1}\\&=
a\frac{M^2}{d^2(1-\lambda)^2}(\lambda
d-\frac{\lambda^{n+1}d}{(n+1)\theta^n}).
\end{split}
\end{equation}
Moreover from Lemma~\ref{lema8'}, we conclude that
\begin{equation}\label{2121}\begin{split}
\frac{a}{\omega_{n-1}}\int_{|y-x_1|\le  d\theta}&|\nabla
u|^2dy\\&\le (n-2)a\frac{d^{n-2}\theta^{n-2}}{1-\theta^{n-2}}
\left(\frac{\gamma\rho^2}{2\beta
n}+\frac{\alpha}{\beta}\max_{|y-x_1|=d}|u(y)-u(x_1)|\right)\\&\le
(n-2)a\frac{d^{n-2}\theta^{n-2}}{1-\theta^{n-2}}
\left(\frac{\gamma\rho^2}{2\beta n}+\frac{2K\alpha}{\beta}\right).
\end{split}
\end{equation}

Where $K:=\max_{|x-x_0|\le r_0} |u(x) - u(x_0)|$.

 Inserting now (\ref{2120}) and (\ref{2121}) in (\ref{2119}) we
obtain
\begin{equation}\label{2122}\begin{split}
&\frac{a}{\omega_{n-1}}\int_{|y-x_1|\le
d\theta}\left(\frac{1}{|y-x_1|^{n-1}}-\frac{|y-x_1|}{d^n\theta^n}\right)|\nabla
u|^2 dy\\&\le a\frac{M^2}{d(1-\lambda)^2}(\lambda
-\frac{\lambda^{n+1}}{(n+1)\theta^n})\\&+(n-2)a\lambda
\left(\frac{1}{d\lambda^{n}}-\frac{1}{d\theta^{n}}\right)\frac{\theta^{n-2}}{1-\theta^{n-2}}
\left(\frac{\gamma\rho^2}{2\beta n}+\frac{2K\alpha}{\beta}\right).
\end{split}
\end{equation}
Combining (\ref{2117}) (for $\theta<\delta_u(\varepsilon)/r_0$),
(\ref{2118}) and (\ref{2122}) we conclude from (\ref{2116}) that the
following inequality holds:
\begin{equation}\label{2123}
\begin{split} \frac{M}{d}&\le
\frac{\min\{n\varepsilon,\gamma_n
\}}{d\theta}+b\frac{n}{n+1}d\theta+a\frac{M^2}{d(1-\lambda)^2}(\lambda
-\frac{\lambda^{n+1}}{(n+1)\theta^n})\\&+
\frac{(n-2)a}{d(1-\theta^{n-2})}\frac{\lambda}{\theta^2} \left(
\left( \frac{\theta}{\lambda} \right)^n-1 \right)
\left(\frac{\gamma\rho^2}{2\beta n}+\frac{2K\alpha}{\beta}\right).
\end{split}
\end{equation}
Myltiplying by $d$ we get:
\begin{equation}\label{2124}
\begin{split}
 M&\le
\frac{\min\{n\varepsilon,\gamma_n
\}}{\theta}+b\frac{n}{n+1}d^2\theta+a\frac{M^2}{(1-\lambda)^2}(\lambda
-\frac{\lambda^{n+1}}{(n+1)\theta^n})\\&+
\frac{(n-2)a}{(1-\theta^{n-2})}\frac{\lambda}{\theta^2} \left(
\left( \frac{\theta}{\lambda} \right)^n-1 \right)
\left(\frac{\gamma\rho^2}{2\beta n}+\frac{2K\alpha}{\beta}\right).
\end{split}
\end{equation}
Remember that $\lambda$ and $\theta$ are arbitrary numbers
satisfying $0<\lambda<\theta<1$. The inequality (\ref{2124}) can be
written in the form
\begin{equation}\label{2125}
AM^2-M+B\ge 0
\end{equation}
where $$A=a\frac{1}{(1-\lambda)^2}(\lambda
-\frac{\lambda^{n+1}}{(n+1)\theta^n})$$ and

\[\begin{split}B&=\frac{\min\{n\varepsilon,\gamma_n \}}{\theta}+b\frac{n}{n+1}d^2\theta\\&+
\frac{(n-2)a}{(1-\theta^{n-2})}\frac{\lambda}{\theta^2} \left(
\left( \frac{\theta}{\lambda} \right)^n-1 \right)
\left(\frac{\gamma\theta^2 d^2}{2\beta
n}+\frac{2K\alpha}{\beta}\right).\end{split}\]  Taking $\lambda=\sin
\theta$ we obtain that $$\lim_{\theta\to 0}4AB\le
\min\{\frac{4a\varepsilon n^2}{n+1},\frac{4a n\gamma_n}{n+1}\}.$$
Hence
$$4AB<1\text{ for }\varepsilon=\frac{n+1}{4an^2+1},$$ whenever
$\theta\le\theta_0$, where $\theta_0$ is small enough. Observe that
in the case
\begin{equation}\label{ieq}\frac{4a  n\gamma_n}{n+1}<1,\end{equation} $\theta_0$ can be chosen
independently of $\varepsilon$ i.e. independently of $u$. The
inequality (\ref{2125}) is equivalent with
\begin{equation}\label{august30}M\le\frac{1-\sqrt{1-4AB}}{2A}=M^-(\theta)\vee M\ge
\frac{1+\sqrt{1-4AB}}{2A}=M^+(\theta)\text{ for
$\theta\le\theta_0$}.\end{equation}

From \eqref{august30} it follows that only one of the following
three cases occur:
\begin{enumerate}
\item $M\le M^-(\theta)$, for $\theta\in (0,\theta_0)$;

\item $M\ge M^+(\theta)$, for $\theta\in (0,\theta_0)$;

\item there exist $\theta_1$, $\theta_2\in (0,\theta_0)$ (say
$\theta_1<\theta_2$), such that $M<  M^-(\theta_1)$ and $M>
M^+(\theta_2)> M^-(\theta_2)$.
\end{enumerate}
As $\lim_{\theta\to 0}\frac{1+\sqrt{1-4AB}}{2A}=+\infty$, the case
(2) is not possible. Since $M^+$ and $M^-$ are continuous, the case
(3) implies that there exists $\theta_3\in (\theta_1,\theta_2)$ such
that $M^-(\theta_3)<M<M^+(\theta_3)$. Thus the case (3) is also
excluded.

The conclusion is that only the case (1) is true and henceforth
$$M\le\frac{1-\sqrt{1-4AB}}{2A}=\frac{2B}{1+\sqrt{1-4AB}}=C_2'(K,\theta_0,a,b,\alpha,\beta,\gamma,r_0,n).$$
Since $d<r_0<1$ it follows that $d^2\le r_0$. Therefore
\begin{equation}\label{final} M\le
2B\le C_1(a,b,\alpha,\beta,\gamma,n,u)(K+r_0).
\end{equation}
From $r_0|\nabla u(x_0)|\le M$ it follows the desired inequality.
\end{proof}

The following two lemmas, roughly speaking assert that the boundary
behavior of any solution of the Poisson differential inequality is
approximately the same as the boundary behavior of the set of two
harmonic mappings. They are $n-$ dimensional "generalizations" of
\cite[Lemma~9]{EH} and \cite[Lemma~9']{EH}. Since the proofs in
\cite{EH} only rely on the maximum principle, the proofs of these
lemmas clearly apply to $n
> 2$ as well with very small modifications.
\begin{lemma}\label{lemma9}
Let $u:B^n\to B^m$ be a $C^2$ mapping defined on the unit ball and
satisfying the inequality
\begin{equation}\label{221} |\Delta u|\le a|\nabla u|^2+b,
\end{equation}
where $0<a<\frac 12$ and $0<b<\infty$. Furthermore let $u(x)$ be
continuous for $|x|\le 1$. Then we have for $x\in B^n$ and $t\in
S^{n-1}$ the estimate
\begin{equation}\label{223}
\begin{split}|u(x)-u(t)|&\le
\frac{1-a}{1-2a}|Y(x)-u(t)|\\&+\frac{a}{2(1-2a)}|F(x)-|u(t)|^2|+\frac{b}{2n(1-2a)}(1-|x|^2),
\end{split}
\end{equation}
where \begin{equation}\label{224}
F(x)=\int_{S^{n-1}}P(x,\eta)|u(\eta)|^2d\sigma(\eta),\, |x|<1
\end{equation}
and
\begin{equation}\label{225}
Y(x)=\int_{S^{n-1}}P(x,\eta)u(\eta)d\sigma(\eta),\, |x|<1.
\end{equation}
\end{lemma}

\begin{lemma}\label{lema15}
Let $\chi:\overline{B^n}\to [-1,1]$ be a mapping of the class
$C^2(B^n)\cap C(\overline{B^n})$ satisfying the differential
inequality:
\begin{equation}\label{2153}
|\Delta\chi|\le a|\nabla \chi|^2+b
\end{equation}
where $a$ and $b$ are finite constants. Then we have for $x\in B^n$
and $t\in S^{n-1}$ the estimate
\begin{equation}\label{2255}
|\chi(x)-\chi(t)|\le\frac{e^a}{a}\left[|h^{p}(x)-e^{a\chi(t)}|+|h^{m}(x)-e^{-a\chi(t)}|+\frac{2ab}n
e^a(1-|x|)\right]
\end{equation}
where \begin{equation}\label{2256}
h^{m}(x)=\int_{S^{n-1}}P(x,\eta)e^{-a\chi(\eta)}d\sigma(\eta)
\end{equation}
and
\begin{equation}\label{2257}
h^{p}(x)=\int_{S^{n-1}}P(x,\eta)e^{a\chi(\eta)}d\sigma(\eta).
\end{equation}
\end{lemma}

\section{ A priori estimate for a solution to Poisson differential inequality}
\begin{theorem}\label{theorem2}
Let $u:D\to\overline{B^m}$ be a $C^2$ mapping, satisfying the
differential inequality:
\begin{equation}\label{2132}
|\Delta u|\le a|\nabla u|^2+b
\end{equation}
where $0<a<1$ and $0<b<\infty$. Then there exists a constant
$c_2=c_2(a,b,n,u)$ such that for $x_0\in D$ and $r_0 =
\mathrm{dist}(x_0,\partial D)$ there holds
\begin{equation}\label{2133} |\nabla u(x_0)|\le
c_2\left(1+\frac{\max_{|x-x_0|\le r_0}|u(x)-u(x_0)|}{r_0}\right).
\end{equation}
If in addition $a\le C_n$ then $c_2$ can be chosen independent of
$u$ and \eqref{2133} is an {\bf a priori estimate}.

\end{theorem}
\begin{proof}
Let us consider the function $G(u)=|u|^2$ and $\phi(x)=G(u(x))$.
Evidently we have
\begin{equation}\label{2135}
|\nabla G(u)|=|2u|\le 2\,\,\,\text{if}\,\,\, |u|\le 1
\end{equation}
and
\begin{equation}\label{2136}\begin{split}
\Delta\phi&=\sum_{i=1}^mD^2G(u)(\nabla u(x)e_i,\nabla
u(x)e_i)+\left<\nabla G(x),\Delta u(x)\right>\\&=2|\nabla
u|^2+2\left<u, \Delta u\right>. \end{split}
\end{equation}
From (\ref{2132})  we conclude
\begin{equation}\label{2137}\Delta\phi\ge 2(1-a)|\nabla u|^2-2b.
\end{equation}
The conditions of Lemma~\ref{lema8} are therefore satisfied by
taking $\alpha=2$, $\beta=2(1-a)$ and $\gamma=2b$. (\ref{2133})
follows with $c_2(a,b,n,u)=c_1(a,b,\alpha,\beta,\gamma,n,u)$.
\end{proof}

\begin{theorem}\label{themain}
Let $u:\overline{B^n}\to \overline{B^m}$ be continuous in
$\overline{B^n}$, $u|_{B^n}\in C^2$, $u|_{S^{n-1}}\in C^{1,\alpha}$
and satisfy the inequalities

\begin{equation}\label{2240}
|\Delta u|\le a|\nabla u|^2+b, \,\,\, x\in B^n,
\end{equation}
\begin{equation}\label{2241}
|u|_{S^{n-1}}|_{1,\alpha}\le K,
\end{equation}
where $0<a<1/2$ and $0<b, K<\infty$.  Then there exists a fixed
positive number $c_4(a,b,n,K,u)$ such that
\begin{equation}\label{2243} |\nabla u(x)|\le
c_4(a,b,n,K,u),\,\,\,x\in B^n.
\end{equation}
If in addition $a\le C_n$ then $c_4(a,b,n,K,u)$ can be chosen
independently of $u$ and \eqref{2243} is an {\bf a priori estimate}.
\end{theorem}

\begin{proof}
Let $x_0=rt\in B^m$, $t\in S^{m-1}$. From Theorem~\ref{theorem2} we
conclude that the inequality \begin{equation}\label{2244} |\nabla
u(x_0)|\le c_2(a,b,n,u)\left(1+\frac{\max_{|x-x_0|\le
1-r}|u(x)-u(x_0)|}{1-r}\right)
\end{equation}
holds.

We shall estimate the quantity $$\displaystyle Q=\max_{|x-x_0|\le
1-r}|u(x)-u(x_0)|.$$ First of all we have $$|u(x)-u(x_0)|\le
|u(x)-u(t)|+|u(x_0)-u(t)|\,\,\text{ for\,\, $|x|<1$}.$$ Applying now
Lemma~\ref{lemma9} we obtain
\begin{equation}\label{split}\begin{split}
|u(x)-u(x_0)|&\le
\frac{1-a}{1-2a}(|Y(x)-u(t)|+|Y(x_0)-u(t)|)\\&+\frac{a}{2(1-2a)}(|F(x)-|u(t)|^2|+|F(x_0)-|u(t)|^2|)\\&+
\frac{b}{2n(1-2a)}[(1-|x|^2)+(1-|x_0|^2)],
\end{split}
\end{equation}
where the harmonic functions $Y$ and $F$ are defined by (\ref{224})
and (\ref{225}). To continue, we use the following result due to
Gilbarg and H\"ormander see \cite[Theorem~6.1 and Lemma~2.1]{GH},
\begin{proposition}\label{hermander} The Dirichlet
problem $\Delta u=f$ in $\Omega$, $u=u_0$ on $\partial \Omega\in
C^1$ has a unique solution $u\in C^{1,\alpha}$, for every $f\in
C^{0,\alpha}$, and $u_0\in C^{1,\alpha}$, and we have
\begin{equation}\label{62}
|u|_{1,\alpha}\le C(|u_0|_{1,\alpha,
\partial\Omega}+|f|_{0,\alpha})
\end{equation}
where $C$ is a constant.\end{proposition} Applying (\ref{62}) on
harmonic functions $Y$ and $F$, according to (\ref{nninequ}),
(\ref{liinequ}) and (\ref{liinequ1}), we first have
\begin{equation}\label{fe}
|Y(x)-u(t)|+|Y(x_0)-u(t)|\le 2C K(1-r)
\end{equation}
and
\begin{equation}\label{feq}
|F(x)-|u(t)|^2|+|F(x_0)-|u(t)|^2|\le 4CK(1-r).
\end{equation}
Combining (\ref{split}), (\ref{fe}) and (\ref{feq}) we obtain
\begin{equation}\label{split1}\begin{split}
|u(x)-u(x_0)|&\le \left[2CK\frac{1-a}{1-2a}+\frac{4CKa}{2(1-2a)}+
\frac{b}{n(1-2a)}\right](1-r).
\end{split}
\end{equation}
Thus for $$c_3(a,b,K,n)=2CK\frac{1-a}{1-2a}+\frac{4CKa}{2(1-2a)}+
\frac{b}{n(1-2a)}$$ we have $$Q\le c_3(a,b,K,n)(1-r).$$ Inserting
this into (\ref{2244}) we obtain $$|\nabla u(x_0)|\le
c_2(a,b,n,u)(1+c_3(a,b,K,n))=c_4(a,b,K,n,u).$$ Since $x_0$ is
arbitrary point of the unit ball the inequality (\ref{2243}) is
established.
\end{proof}

Whether
 Theorem~\ref{themain} holds replacing the condition $0<a<1/2$ by
$0<a<\infty$, is not known by the author.
However adding the condition of quasiregularity we obtain the
following extension of Theorem~\ref{themain}.
\begin{corollary}
Assume that  $u:\overline{B^n}\to \Bbb R^n$ is a $K-$quasiregular,
twice differentiable mapping, continuous on $\overline{B^n}$, and
$u|_{S^{n-1}}\in C^{1,\alpha}$. If in addition it satisfies the
differential inequality
\begin{equation}\label{qcl}
|\Delta u|\le a|\nabla u|^2+b\text{ for some constants $a,b>0$}
\end{equation}
then $$|\nabla u|\le C_{8}(K,a,b,u).$$
\end{corollary}
\begin{proof}
From (\ref{quasin}) we obtain for $i=1,\dots n$
\begin{equation}\frac{1}{K}\le \frac{l(\nabla u)^n}{J_u}\le
\frac{|\nabla u_i|_2^n}{J_u}\le \frac{|\nabla u|^n}{J_u}\le K
\end{equation}
and hence
\begin{equation}|\nabla u|\le K^{4/n}|\nabla u_i|_2.
\end{equation}
Thus for every $i=1,\dots,n$ \begin{equation} |\Delta u_i|\le
aK^{2/n}|\nabla u_i|_2^2+b.
\end{equation}
The conclusion follows according to Theorem~\ref{themain1}.
\end{proof}

In the rest of the paper we will prove an analogous result for
arbitrary $a$ and $b$. The only restriction is $u$ being a real
function, i.e. $m=1$.

\begin{theorem}\label{theorem2'}
Let $B^n(x_0,r_0)\subset D\subset B^n$ and let
$\chi:D\subset{B^n}\to [-1,1]$ be a mapping of the class $C^2(D)$
satisfying the differential inequality:
\begin{equation}\label{2148}
|\Delta\chi|\le a|\nabla \chi|^2+b
\end{equation}
where $a$ and $b$ are finite constants. Then we have the estimate
\begin{equation}\label{2250}
|\nabla \chi(x_0)|\le c_5(a,b,n,\chi)\left(1+\frac{\max_{|x-x_0|\le
r_0}|\chi(x)-\chi(x_0)|}{r_0}\right).
\end{equation}
If $a\le C_n$ then $c_5(a,b,n,\chi)$ can be chosen independently of
$\chi$ and \eqref{2250} is an {\bf a priori estimate}.
\end{theorem}
\begin{proof}
Let us consider a twice differentiable function $\phi(t), $ $-1\le
t\le 1$ and $\varphi(x)=\phi(\chi(x))$. The function $\varphi$
satisfies the differential equation
\begin{equation}\label{2151}\Delta\varphi=\phi''|\nabla\chi|^2+\phi'\Delta\chi.\end{equation} Using
(\ref{2148}) we obtain
\begin{equation}\label{2151'}
\Delta\varphi\ge (\phi''-a|\phi'|)|\nabla\chi|^2-b|\phi'|.
\end{equation}
Taking $\phi(t)=e^{2at}$ we obtain
\begin{equation}\label{2152}
\Delta \varphi\ge 2a^2e^{-2a}|\nabla \chi|^2-2abe^{2a}.
\end{equation}
The conditions of Lemma~\ref{lema8} are therefore satisfied by
taking $\alpha=2ae^{2a}$, $\beta=2a^2e^{-2a}$, and
$\gamma=2abe^{2a}$. Hence we conclude that (\ref{2250}) holds for
$c_5(a,b,n,\chi)=c_1(a,b,\alpha,\beta,\gamma,n,\chi)$.
\end{proof}

\begin{theorem}\label{themain1}
Let $\chi:\overline{B^n}\to \Bbb R$ be continuous in
$\overline{B^n}$, $\chi|_{B^n}\in C^{2}$, $\chi|_{S^{n-1}}\in
C^{1,\alpha}$ and satisfy the inequalities

\begin{equation}\label{2240'}
|\Delta \chi|\le a|\nabla \chi|^2+b, \,\,\, x\in B^n,
\end{equation}
\begin{equation}\label{2241'}
|\chi|_{S^{n-1}}|_{1,\alpha}\le K
\end{equation}
where $0<a, b, K$.  Then there exists a fixed positive number
$c_6=c_6(a,b,K,n,\chi)$, which do not depends on $\chi$ for  $a\le
C_n|/|\chi|_\infty$ such that
\begin{equation}\label{2243'} |\nabla \chi(x)|\le
c_6,\,\,\,x\in B^n.
\end{equation}
\end{theorem}
\begin{remark}The condition $|\chi|_{S^{n-1}}|_{1,\alpha}\le K$ of
Theorem~\ref{themain1} is the best possible, i.e. we cannot replace
it by  $|\chi|_{S^{n-1}}|_{1}\le K$. For example O. Martio in
\cite{OM} gave an example of a harmonic diffeomorphism $w=P[f]$, of
the unit disk onto itself such that $f\in C^1(S^1)$ and  $\nabla w$
is unbounded. This example can be easily modified for the space. For
example we can simply take $u(x_1,x_2,\dots,x_n)=P[f](x_1,x_2)$.
Then $u|_{S^{n-1}}\in C^1$ but $\nabla u$ is not bounded.
\end{remark}
\begin{proof}
The proof follows the same lines as the proof of
Theorem~\ref{themain}. The only difference is applying
Theorem~\ref{theorem2'} instead of Theorem~\ref{theorem2} and
Lemma~\ref{lema15} instead of Lemma~\ref{lemma9} to the function
$\chi_0=\chi(x)/M$, where  $M=\max\{|\chi(x)|: x\in
\overline{B^{n}}\}$.

Let $x_0=rt\in B^m$, $t\in S^{m-1}$. From Theorem~\ref{theorem2'} we
conclude that the inequality
\begin{equation}\label{2244'}
|\nabla \chi_0(x_0)|\le
c_5(a,b,n,\chi)\left(1+\frac{\max_{|x-x_0|\le
r_0}|\chi_0(x)-\chi_0(x_0)|}{r_0}\right)
\end{equation}
holds.

We shall estimate the quantity $$\displaystyle Q=\max_{|x-x_0|\le
1-r}|\chi_0(x)-\chi_0(x_0)|.$$ First of all we have
\begin{equation}\label{foal}|\chi_0(x)-\chi_0(x_0)|\le
|\chi_0(x)-\chi_0(t)|+|\chi_0(x_0)-\chi_0(t)|\,\,\text{ for\,\,
$|x|<1$}.\end{equation} Applying now Lemma~\ref{lema15} we obtain
\begin{equation}\label{split8}\begin{split}
|\chi_0(x)-\chi_0(x_0)|&\le\frac{e^a}{a}\bigg[|h^{p}(x)-e^{a\chi_0(t)}|+|h^{p}(x_0)-e^{a\chi_0(t)}|\bigg]
\\&+\frac{e^a}{a}\bigg[|h^{m}(x)-e^{-a\chi_0(t)}|+|h^{m}(x_0)-e^{-a\chi_0(t)}|\bigg]\\&+\frac{2ab}{n}e^a(1-|x|+1-|x_0|),
\end{split}
\end{equation}
where the harmonic functions $h^p$ and $h^m$ are defined by
(\ref{2256}) and (\ref{2257}). Applying (\ref{62}) on harmonic
functions $h^p$ and $h^m$, according to (\ref{nninequ}),
(\ref{liinequ}) and (\ref{liinequ2}) for $\tau_1=e^{at}$ and
$\tau_2=e^{-at}$, we first have
\begin{equation}\label{fe1}
|h^{p}(x)-e^{a\chi_0(t)}|+|h^{p}(x_0)-e^{a\chi_0(t)}|\le 2ae^aC
K(1-r)
\end{equation}
and
\begin{equation}\label{feq1}
|h^{m}(x)-e^{-a\chi_0(t)}|+|h^{m}(x_0)-e^{-a\chi_0(t)}|\le 2ae^aC
K(1-r).
\end{equation}
Combining (\ref{split8})-(\ref{feq1}) we obtain
\begin{equation}\label{split10}\begin{split}
|\chi_0(x)-\chi_0(x_0)|&\le \left[4CKe^{2a}+
\frac{4a^2b}{n}e^a\right](1-r).
\end{split}
\end{equation}
Thus for $$c_7(a,b,K,n)=4CKe^{2a}+ \frac{4a^2b}{n}e^a$$ we have
$$Q\le C_7(a,b,K)(1-r).$$ Inserting this into (\ref{2244'}) we
obtain $$|\nabla \chi_0(x_0)|\le
c_5(a,b,n,\chi)(1+c_7(a,b,K,n))=c'_6(a,b,K,n,\chi).$$ Since $x_0$ is
an arbitrary point of the unit ball the inequality (\ref{2243'}) is
valid for $c_{6}(a,b,K,n,\chi)=c'_{6}(a,b,K,n,\chi)\cdot M$.
\end{proof}

\section{Applications-The proof of Theorem~C}
\subsection{Bounded curvature and the distance function}
Let $\Omega$ be a domain in $\Bbb R^n$ having a non-empty boundary
$\partial \Omega$. The distance function is defined by
\begin{equation}\label{dist}
d(x)=\mathrm{dist}\,(x,\partial\Omega).
\end{equation}
The function $d$ is uniformly Lipschitz continuous and there holds
the inequality \begin{equation}\label{distine} |d(x)-d(y)|\le|x-y|.
\end{equation}
Now let $\partial \Omega\in C^2$. For $y\in \partial \Omega$, let
$\mathbf{\nu}(y)$ and $T_y$ denote respectively the unit inner
normal to $\partial\Omega$ at $y$ and the tangent hyperplane to
$\partial\Omega$ at $y$.

The curvature of $\partial\Omega$ at a fixed point $y_0\in
\partial\Omega$ is determined as follows. By the rotation of
coordinates we can assume that $x_n$ coordinate axis lies in the
direction $\mathbf{\nu}(y_0)$. In some neighborhood $\mathcal
N(y_0)$ of $y_0$, $\partial\Omega$ is given by $x_n=\varphi(x')$,
where $x'=(x_1,\dots,x_{n-1})$, $\varphi\in C^2(T(y_0)\cap \mathcal
N(y_0))$ and $\nabla\varphi(y_0')=0$. The curvature of
$\partial\Omega$ at $y_0$ is then described by orthogonal invariants
of the Hessian matrix $D^2\varphi$ evaluated at $y_0$. The
eigenvalues of $D^2\varphi(y_0')$, $\kappa_1$, $\dots$,
$\kappa_{n-1}$ are called the principal curvatures of
$\partial\Omega$ at $y_0$ and the corresponding eigenvectors are
called the principal directions of $\partial\Omega$ at $y_0$. The
mean curvature of $\partial\Omega$ at $y_0$ is given by
\begin{equation}\label{meanc}
H(y_0)=\frac{1}{n-1}\sum_{i=1}^{n-1}\kappa_i=\Delta\varphi(y_0').
\end{equation}
By a further rotation of coordinates we can assume that the
$x_1,\dots,x_{n-1}$ axes lie along principal directions
corresponding to $\kappa_1$, $\dots$, $\kappa_{n-1}$ at $y_0$. The
Hessian matrix $D^2\varphi(y_0')$ with respect to the principal
coordinate system at $y_0$ described above is given by
$$D^2\varphi(y_0')=\mathrm{diag}(\kappa_1,\dots,\kappa_{n-1}).$$
\begin{proposition}\cite{gt}
Let $\Omega$ be bounded domain of class $C^k$ for $k\ge 2$. Then
there exists a positive constant $\mu$ depending on $\Omega$ such
that $d\in C^k(\Gamma_\mu)$, where
$\Gamma_\mu=\{x\in\overline{\Omega}:d(x)<\mu\}$ and for $x\in
\Gamma_\mu$ there exists $y(x)\in \partial\Omega$ such that
\begin{equation}\label{distnorm} \nabla d(x)=\mathbf{\nu}(y(x)).
\end{equation}
\end{proposition}
\begin{proposition}\cite{gt}
Let $\Omega$ be of class $C^k$ for $k\ge 2$. Let $x_0\in\Gamma_\mu$,
$y_0\in\partial\Omega$ be such that $|x_0-y_0|=d(x_0)$. Then in
terms of a principal coordinate system at $y_0$, we have
\begin{equation}\label{diag}
D^2d(x_0)=\mathrm{diag}(\frac{-\kappa_1}{1-\kappa_1d},\dots,\frac{-\kappa_{n-1}}{1-\kappa_{n-1}d},0).
\end{equation}
\end{proposition}
\begin{lemma} Let $\partial\Omega\in C^2$. For $x\in \Gamma_\mu$ and $y(x)\in \partial\Omega$ there holds the equation
\begin{equation}\label{deltad} \Delta
d(x)= \sum_{i=1}^{n-1}\frac{-\kappa_i(y(x))}{1-\kappa_i(y(x))d}.
\end{equation}
If for some $x_0\in \Omega$, the mean curvature of $y_0=y(x_0)\in
\partial \Omega$, is positive: then
$-d(x)$ is subharmonic in some neighborhood of $y_0$. In particular
if  $\Omega$ is convex then the function $\Gamma_\mu\ni x\mapsto
-d(x)$ is subharmonic.
\end{lemma}
\begin{proof}
The equation (\ref{deltad}) follows from (\ref{diag}) and
(\ref{hess}) (it is a special case of the relation (\ref{ende})
taking $u(x)=x$). If $\Omega$ is convex then for every $i$
$\kappa_i\ge 0$. Hence $\Delta(-d(x))\ge 0$ and thus $-d(x)$ is
subharmonic.
\end{proof}

\begin{lemma}
Let $u:\Omega\to \Omega'$ be a $K$ q.r. and $\chi=-d(u(x))$. Then
\begin{equation}\label{dqh} |\nabla\chi|\le|\nabla u|\le K^{2/n}|\nabla
\chi|
\end{equation}
in $u^{-1}(\Gamma_\mu)$ for $\mu>0$ such that $1/\mu>\kappa_0=
\max\{|\kappa_i(x)|:x\in\partial\Omega,i=1,\dots,n-1\}$.
\end{lemma}
\begin{proof}
Observe first that $\nabla d$ is a unit vector. From $\nabla
\chi=-\nabla d\cdot \nabla u$ it follows that
$$|\nabla\chi|\le|\nabla d| |\nabla u|=|\nabla u|.$$ To continue we
need the following observation. For a non-singular matrix $A$ we
have
\begin{equation}\begin{split}\inf_{|x|=1}
|Ax|^2&=\inf_{|x|=1}\left<Ax,Ax\right>=\inf_{|x|=1}\left<A^t A
x,x\right>\\&=\inf\{\lambda: \exists x\neq 0, A^t A x=\lambda
x\}\\&=\inf\{\lambda: \exists x\neq 0, A A^t A x=\lambda
Ax\}\\&=\inf\{\lambda: \exists y\neq 0, A A^t y=\lambda
y\}=\inf_{|x|=1} |A^t x|^2.\end{split}\end{equation} Next we have
$\nabla \chi=-(\nabla u)^t\cdot \nabla d$ and therefore for $x\in
u^{-1}(\Gamma_\mu)$, we obtain $$|\nabla \chi|\ge
\inf_{|e|=1}|(\nabla u)^t \,e|=\inf_{|e|=1}|\nabla u \,e|=l(u)\ge
K^{-2/n}|\nabla u|.$$

The proof of (\ref{dqh}) is completed.

\end{proof}

\begin{lemma}\label{lemo} Let $D$ and $\Omega\subset\Bbb R^n$ be open domains, and
$\partial\Omega$ be a $C^2$ hypersurface homeomorphic to $S^{n-1}$.
Let $u:D\to \Omega$ be a twice differentiable $K$ quasiregular
surjective mapping satisfying the Poisson differential inequality.
Let in addition $\chi(x)=-d(u(x))$. Then there exists a constant
$a_1=C(a,b,K,\Omega)$ such that
$$|\Delta\chi(x)|\le a_1|\nabla \chi(x)|^2+b$$ in
$u^{-1}(\Gamma_\mu)$ for some $\mu>0$ with
$1/\mu>\kappa_0=\max\{|\kappa_i(y)|:y\in\partial\Omega,i=1,\dots,n-1\}$.
If in addition $u$ is harmonic and $\Omega$ is convex then $\chi$
is subharmonic for $x\in u^{-1}(\Gamma_\mu)$.
\end{lemma}
\begin{proof}
Let $y\in\partial\Omega$. By the considerations taken in the begin
of this section we can choose an orthogonal transformation $O_y$ so
that the vectors $O_y(e_i)$, $i=1,\dots,n-1$ make the principal
coordinate system in the tangent hyperplane $T_{y}$ of
$\partial\Omega$, that determine the principal curvatures of
$\partial\Omega$ and $O_y(e_n)=\mathrm{\nu}(y)$. Let $x_0\in B^n$.
Choose $y_0\in \partial \Omega$ so that
$d(u(x_0))=\operatorname{dist}\,(u(x_0),y_0)$. Take
$\tilde\partial\Omega:=O_{y_0}\partial\Omega$. Let $\tilde d$ be the
distance function with respect to $\tilde\partial\Omega$. Then
$d(u)=\tilde d(O_{y_0}(u))$ and $\chi(x)=-\tilde d(O_{y_0}(u(x)))$.
Thus
\begin{equation}\label{hess}\begin{split}
\Delta\chi(x)&=-\sum_{i=1}^nD^2\tilde d(O_{y_0}(u(x)))(O_{y_0}\nabla
u(x)e_i,O_{y_0}\nabla u(x)e_i)\\&-\left<\nabla d(u(x)),\Delta
u(x)\right>.
\end{split}
\end{equation}

Next we have \begin{equation}\label{dlta}|\left<\nabla
d(u(x)),\Delta u(x)\right>|\le |\Delta u|\le a|\nabla
u|^2+b.\end{equation} Applying (\ref{diag})
\begin{equation}\label{ende}\begin{split}
\sum_{i=1}^{n}D^2&\tilde d(O_{y_0}(u(x_0)))(O_{y_0}(\nabla
u(x_0))e_i,O_{y_0}(\nabla u(x_0))e_i)\\&=\sum_{i=1}^n
\sum_{j,k=1}^nD_{j,k}\tilde d(O_{y_0}(u(x_0))) \,D_i
(O_{y_0}u)_j(x_0)\cdot D_i (O_{y_0}u)_k(x_0)\\&=\sum_{j,k=1}^n
D_{j,k}\tilde d (O_{y_0}(u(x_0)))\left<(O_{y_0} \nabla
u(x_0))^te_j,(O_{y_0} \nabla u(x_0))^te_k\right>\\&=
\sum_{i=1}^{n-1}\frac{-\tilde\kappa_i}{1-\tilde\kappa_i\tilde
d}|(O_{y_0} \nabla u(x_0))^te_i|^2.
\end{split}
\end{equation}
Since the principal curvatures $\tilde \kappa_i=\kappa_i$ are
bounded by $\kappa_0$, combining (\ref{hess}), (\ref{dlta}),
(\ref{ende}) and (\ref{dqh}), and using the relations $$|(O_{y_0}
\nabla u(x_0))^te_i|^2=|(\nabla u(x_0))^tO^t_{y_0}e_i|^2\le |\nabla
u(x_0)|^2,$$ we obtain for $x\in u^{-1}(\Gamma_\mu)$
\begin{equation}
|\Delta\chi|\le K^{4/n}(a+\frac{n\kappa_0}{1-\mu\kappa_0})|\nabla
\chi|^2+b,
\end{equation}
which is the desired inequality.
\end{proof}
A $K$ q.c. self-mapping $f$ of the unit ball $B^n$ need not be
Lipschitz continuous. It is holder continuous  i.e. there hold the
inequality
\begin{equation}|f(x)-f(y)|\le
M_1({n,K})|x-y|^{K^{1/(1-n)}}.
\end{equation}
See \cite{fv} for the details. See \cite{pjj} for the extension of
the Mori's theorem for domains satisfying the quasihyperbolic
boundary conditions as well as for quasiconvex domains.

Under some additional conditions on interior regularity we obtain
that a q.c. mapping is Lipschitz continuous.
\begin{theorem}[The main result]\label{arbitrarydomain}
Let $u:B^n\to \Omega$ be a twice differentiable quasiconformal
mapping of the unit ball onto the bounded domain $\Omega$ with $C^2$
boundary satisfying the Poisson differential inequality. Then
$\nabla u$ is bounded and $u$ is Lipschitz continuous.
\end{theorem}

\begin{proof}
From Lemma~\ref{lemo} $$|\Delta \chi|\le a_1|\nabla
\chi|^2+b_1\text{ for $x\in \Gamma_\mu$}.$$ On the other hand, by a
theorem of Martio and Nyakki (\cite{mar}) $u$ has a continuous
extension to the boundary. Therefore  for every $x\in S^{n-1}$,
$\lim_{y\to x} \chi(y)=\chi(x)=0$. Let $\tilde\chi$ be an $C^2$
extension of the function $\chi|_{x\in u^{-1}(\Gamma_\mu)}$ in $B^n$
(by Whitney theorem it exists \cite{whh}). Let $b_0=\max\{|\Delta
\tilde\chi(x)|:x\in B^n\setminus u^{-1}(\Gamma_{\mu/2})\}$. Then
$$|\Delta \tilde \chi| \le a_1|\nabla \tilde\chi|^2+b_1+b_0.$$ Thus
the conditions of Theorem~\ref{themain1} are satisfied. The
conclusion is that $\nabla \tilde\chi$ is bounded. According to
(\ref{dqh}) $\nabla u$ is bounded in $u^{-1}(\Gamma_\mu)$ and hence
in $B^n$ as well. The conclusion of the theorem now easily follows.
\end{proof}
Let $u=P[f]$. Let
$S=S(r,\theta)=S(r,\varphi,\theta_1,\dots,\theta_{n-2})$,
$\theta\in[0,2\pi]\times[0,\pi]\times\dots\times[0,\pi]$ be the
spherical coordinates and let $T(\theta)=S(1,\theta)$. Let in
addition $x=f(T(\theta))$ and $\mathbf{n}_x$ be the normal on
$\partial\Omega$ defined by the formula
$\mathbf{n}_x=x_{\varphi}\times x_{\theta_1}\times\dots\times
x_{\theta_{n-2}}$. Since $f(S^{n-1})=\partial \Omega$ it follows
that
\begin{equation}\label{normal}
\mathbf{n}_x=|\mathbf n_x| \nu_x=D_x \cdot \nu_x
\end{equation}
where $\nu_x$ is the unit inner normal vector that defines the
tangent hyperplane of $\partial \Omega$ at $x =f(T(\theta))$
$$TP^{n-1}_{x} = \{y: \left<x-y,\mathbf
\nu_x\right>=0\}.$$ Since $\Omega$ is convex it follows that
\begin{equation}\label{convexo} \left<x-y,\mathbf \nu_x\right>\ge 0
\text{ for every $x\in
\partial\Omega$ and $y\in \Omega$}.
\end{equation}
 Let in addition
$u(S(r,\theta))=(y_1,y_1,\dots,y_n)$. Then in these terms we have
the following corollary.
\begin{corollary}
If $u$ is a q.c. harmonic mapping of the unit ball onto a convex
domain $\Omega$ with $C^2$ boundary, then:
\begin{equation}\label{boundjac}J_u\in
L^\infty(B^n),\end{equation}
\begin{equation}\label{jas1}J^b_u(t):=\lim_{r\to 1}J_u(rt)\in L^\infty(S^{n-1}),\end{equation} and there hold the inequality
\begin{equation}\label{jac2}J^b_u(t) \ge
\frac{\mathrm{dist}(u(0),\partial\Omega)^n}{(2K)^{n^2-n}},\end{equation}
where $K$ is the quasiconformality constant.

\end{corollary}

We need the following lemma.
\begin{lemma}\label{prim}Let $A:\Bbb R^n\rightarrow \Bbb R^n$ be a linear operator such that
$A=[a_{ij}]_{{i,j =1,\dots, n}}$. If $A$ is $K$ quasiconformal, then
there hold the following double inequality
\begin{equation}\label{useful}K^{1-n}|A|^{n-1}|x_1\times\dots\times
x_{n-1}|\le |A x_1\times \dots\times Ax_{n-1}|\le
|A|^{n-1}|x_1\times\dots\times x_{n-1}|.
\end{equation}
Here $\times\dots\times $ denotes the vectorial product. Both
inequalities in \eqref{useful} are sharp.
\end{lemma}
The author believes that the Lemma~\ref{prim} is well-known, and its
proof is given in the forthcoming author's paper \cite{kja}.
\begin{proof}
 Since (in view of
Theorem~\ref{arbitrarydomain}) $\nabla u=(D_iu_j)_{i,j=1}^n$ is
bounded, every harmonic mapping $D_iu_j$ is bounded. Therefore there
exists $v_{i,j}\in L^{\infty}(S^{n-1})$ such that
$D_iu_j=P[v_{i,j}]$. Thus $\lim_{r\to 1}D_iu_j(rt)=v_{i,j}(t)$ for
every $i,j$. The relations (\ref{boundjac}) and (\ref{jas1}) are
therefore proved. On the other hand since $y_j=u_j\circ S$,
$j=1,\dots, n$, we have for a.e. $t=S(1,\theta)\in S^{n-1}$ the
relations:
\begin{equation}\label{uno}
\lim_{r \to 1}{y_i}_{\varphi}(r,\theta)={x_i}_{\varphi}(\theta),\
i\in \{1,\dots , n\},
\end{equation}

\begin{equation}\label{due}
\lim_{r \to 1}{y_i}_{\theta_j}(r,\theta)={x_i}_{\theta_j}(\theta),\
i\in \{1,\dots , n\},\, j\in \{1,\dots, n-2\},
\end{equation}

and

\begin{equation}\label{tre}
\lim_{r \to 1}{y_i}_r(r,\theta)=\lim_{r \to
1}\frac{x_i(\theta)-y_i(r,\theta)}{1-r},\ i \in \{1,\dots , n\}.
\end{equation}
From (\ref{uno}), (\ref{due}), (\ref{tre}) and (\ref{green}) we
obtain for a.e. $t=S(1,\theta)\in S^{n-1}$:

\begin{equation*}\label{equstar}
\begin{split}
\lim_{r \to 1}J_{u\circ S}(r,\theta) &=\lim_{r \to 1} \left|
\begin{array}{cccc}
\frac{x_1-y_1}{1-r} & \frac{x_2-y_2}{1-r} & \dots &
\frac{x_n-y_n}{1-r}\\ {x_1}_\varphi &{x_2}_\varphi & \dots &
{x_n}_{\varphi}\\ {x_1}_{\theta_1} &{x_2}_{\theta_1} & \dots &
{x_n}_{\theta_1}\\ \hdotsfor[2]{4} \\ {x_1}_{\theta_{n-2}}
&{x_2}_{\theta_{n-2}} & \dots & {x_n}_{\theta_{n-2}}
\end{array}
\right|
\\
&= \lim_{r \to 1}\int_{S^{n-1}}\frac{1+r}{|\eta-x|^n} \left|
\begin{array}{ccc}
{x_1-f_1(\eta)} & \dots & {x_n-f_n(\eta)}\\ {x_1}_{\varphi} & \dots
& {x_n}_{\varphi}\\ {x_1}_{\theta_1} & \dots & {x_n}_{\theta_1}\\
\hdotsfor[2]{3} \\ {x_1}_{\theta_{n-2}} & \dots &
{x_n}_{\theta_{n-2}}
\end{array}\right| d\sigma(\eta)\\&=\lim_{r \to 1}\int_{S^{n-1}}\frac{1+r}{|\eta-S(r,\theta)|^n}
\left<f(T(\theta))-f(\eta),{\mathbf n}_{f\circ T
}(T(\theta))\right>d\sigma(\eta).
\end{split}
\end{equation*}

Using (\ref{convexo}) and the inequality $$\lim_{r\to 1}
\frac{1+r}{|\eta-S(r,\theta)|^n}\ge \frac{1}{2^{n-1}}$$ we obtain
\begin{equation*}
\begin{split}
\lim_{r \to 1}J_{u\circ S}(r,\theta)&\ge
\frac{D_x(\theta)}{2^{n-1}}\int_{S^{n-1}}\left<f(T(\theta))-f(\eta),\nu_x\right>d\sigma(\eta)
\\&=\frac{D_x(\theta)}{2^{n-1}}
\left(\left<f(T(\theta)),\nu_x\right>-\left<u(0),\nu_x\right>\right)\\
&=\frac{D_x(\theta)}{2^{n-1}}\left<f(T(\theta))-u(0),\nu_x\right>
\\
&=\frac{D_x(\theta)}{2^{n-1}}\mathrm{dist}\left(TP^{n-1}_{f(S(1,\theta))},u(0)\right)\\
&\geq\frac{D_x(\theta)}{2^{n-1}}\mathrm{dist}(u(0),\partial\Omega).
\end{split}
\end{equation*}

Thus for a.e. $t=S(1,\theta)\in S^{n-1}$, we have
\begin{equation}\label{jacobi}J^b_u(S(1,\theta))=\frac{J_{u\circ
S}(\theta)}{D_T(\theta)}\ge
\frac{D_x(\theta)}{D_T(\theta)}\frac{\mathrm{dist}(u(0),\partial\Omega)}{2^{n-1}}.\end{equation}

From the left side of \eqref{useful}, using the inequality $$|\nabla
u|^n\ge J_u(t)$$ we obtain $$\frac{D_x(\theta)}{D_T(\theta)} \ge
K^{1-n}|\nabla u(t)|^{n-1}\ge K^{1-n} J_u(t)^{(n-1)/n}.$$ Combining
the last inequality and \eqref{jacobi} we obtain (\ref{jac2}).
\end{proof}
\subsection{An open problem}  It remains an
open problem whether every q.c. harmonic mapping of the unit ball
onto a domain with $C^2$ boundary is bi-Lipschitz continuous. This
question has affirmative answer for the plane case (see
\cite{kalajan}).

\

\end{document}